\documentclass[12pt]{article}
\usepackage{latexsym}
\font\smallit=cmti10

\begin{document} 
\begin{center}
{\bf \Large A novel operation associated with Gauss' \\
\vspace{0.5cm}
arithmetic-geometric means}
\vskip 20pt
{\bf Shinji Tanimoto}\\
{\smallit Department of Mathematics, Kochi Joshi University, Kochi 780-8515, Japan}\\
{\tt tanimoto@cc.kochi-wu.ac.jp}\\ 
\end{center}
\begin{abstract}
The arithmetic mean is the mean for addition and the geometric mean is that for multiplication.
Then what kind of binary operation is associated with the arithmetic-geometric mean (AGM) due to C. F. Gauss?
If it is possible to construct an arithmetic operation such that AGM is the mean for this operation, 
it can be regarded as
an intermediate operation between addition and multiplication in view of the meaning of AGM.
In this paper such an operation is introduced and several of its algebraic properties are proved.
\end{abstract}
\begin{center}
{\bf 1. Introduction}
\end{center}

This paper is an enlarged English version of [4] that was published in Japanese ten years ago. 
The objective is to add several remarks to it and to be made available by a wider readership. \\
\indent
The arithmetic mean is the mean for addition and 
the geometric mean is that for multiplication.
Then what kind of binary operation is associated with the arithmetic-geometric mean (AGM) due to C. F. Gauss?
The AGM of two positive real numbers $x$ and $y$, denoted by ${\rm agm} (x, y)$,
is defined as follows (see [2]). Let the sequences be given by
\begin{eqnarray*}
x_0 = x, ~ y_0 = y, ~ x_n = (x_{n-1} + y_{n-1})/2, ~ y_n = \sqrt{x_{n-1}y_{n-1}} ~~ (n \ge 1).
\end{eqnarray*} 
Then it is defined as their common limits:
\begin{eqnarray}
{\rm agm} (x, y) = \lim_{n\to\infty} x_n = \lim_{n\to\infty} y_n. 
\end{eqnarray}
We propose a new binary operation on the set of all positive real numbers and show that the mean with respect to
the operation gives rise to AGM. That is, if we denote this binary operation by $\star$, the relation is expressed by
\begin{eqnarray}
    x \star y = {\rm agm} (x, y) \star {\rm agm} (x, y).
\end{eqnarray} 
In view of the meaning of AGM, it follows from (2) that the binary operation can be regarded
as an intermediate operation between addition and multiplication. 
Since the equality
\[
{\rm agm} (\lambda x, \lambda y) = \lambda{\rm agm} (x, y)
\] always holds for all $\lambda > 0$, it suffices to treat 
the form ${\rm agm} (1, x)$ for $x >0$. \\
\indent
In order to do this we need one of Jacobi's theta functions (see [1, 3]), simply denoted by $\theta$:
\[
    \theta(q) = \sum_{n = - \infty}^{+ \infty} q^{n^2} = 1 + 2\sum_{n = 1}^{\infty} q^{n^2}
\]
for $-1 <q<1$. 
According to [1], the function $\theta(q)$ increases monotonically in this interval such that
\begin{eqnarray}
\lim_{q\to -1} \theta(q) = 0, ~~ \theta(0) = 1, ~~ \lim_{q\to 1} \theta(q) = + \infty,
\end{eqnarray} 
and satisfies the equality 
\begin{eqnarray}
    {\rm agm} (1, \theta^2(-q)/ \theta^2(q)) =  1/\theta^2(q) ~~ {\rm or}~~
     {\rm agm} (\theta^2(q), \theta^2(-q)) = 1,
\end{eqnarray} 
for all $q~(-1 <q<1)$, where $\theta^2(q)$ denotes the square of $\theta(q)$. 
This equality was found by Gauss himself. \\

\begin{center}
{\bf 2. A binary operation and its properties}
\end{center}

Now we define a new binary operation on the set of all positive real numbers. After that we will
proceed to prove fundamental properties of the operation. Note that,
for positive numbers $x$ and $y$, we have ${\rm agm}(x,y) > 0$ and that there exists 
a unique $q~(-1 <q<1)$ such that $1 / {\rm agm}(x,y) = \theta^2(q)$ from (3).\\
\\
{\bf Definition.} {\it 
For any two positive numbers $x$ and $y,$ choose a unique $q~(-1 <q<1)$ such that 
$1 / {\rm agm}(x,y) = \theta^2(q)$. Then define
\begin{eqnarray}
     x \star y =  \theta^2(-q)/\theta^2(q).
\end{eqnarray} 
}
\\
\indent
This definition together with (4) leads us to 
\begin{eqnarray}
{\rm agm} (1, x \star y) = {\rm agm}(x, y), 
\end{eqnarray}
since ${\rm agm}(x, y) = 1/\theta^2(q)$ holds. It is evident that
$x \star y > 0$ and $x \star y = y \star x$. 
We state fundamental properties of this operation in the following theorem, 
where variables $x, y$ are positive real numbers.\\
\\
\noindent
{\bf Theorem 1.} {\it The operation $\star$ defined above satisfies the following properties.}\\
~{\bf (A)} $1 \star x = x$ {\rm for all} $x$.  {\it Hence {\rm 1} is the unit element of the operation}. \\
~{\bf (B)} $x \star x = y \star y$ {\it implies} $x = y$. \\
~{\bf (C)} $x \star y = {\rm agm} (x, y) \star {\rm agm} (x, y)$. {\it Thus the mean with respect to
the operation is the AGM.} \\
\\
\noindent
{\bf Proof.} \\
~(A)~
Since $\theta^2(-q)/\theta^2(q)$ decreases monotonically in the interval $- 1 < q < 1$ and
\[
\lim_{q\to -1} \theta^2(-q)/\theta^2(q) = \infty,~~ \lim_{q\to 1}\theta^2(-q)/\theta^2(q) = 0,
\]
it is always possible to find $q$ such that $x= \theta^2(-q)/\theta^2(q)$.
Next let $q' (-1 <q'<1)$ be such that $\theta^2(q') = 1/ {\rm agm} (1, x)$. Then by definition we have
$1 \star x = \theta^2(-q')/\theta^2(q')$. On the other hand, it follows from (4) that
\[
{\rm agm} (1, x) = {\rm agm} (1, \theta^2(-q)/ \theta^2(q)) =  1/\theta^2(q). 
\]
Hence we get 
$ \theta^2(q) = \theta^2(q')$. This yields $q' = q$ and thus $1 \star x = x$. \\
~(B)
Let $x^{-1} = \theta^2(q)$ and $y^{-1} = \theta^2(q')$. Since ${\rm agm}(x, x)=x$ and ${\rm agm}(y, y)=y$, 
by definition we have $x \star x = \theta^2(-q)/\theta^2(q)$ and $y \star y = \theta^2(-q')/\theta^2(q')$.
If $x \star x = y \star y$ holds, then $q = q'$ and hence $x = y$. \\
~(C)
Letting $1 / {\rm agm}(x,y) = \theta^2(q)$, we have $x \star y =  \theta^2(-q)/\theta^2(q)$.
Since 
\[
{\rm agm}(x, y)= {\rm agm}({\rm agm}(x, y), {\rm agm}(x, y)),
\]
 it is obvious that
${\rm agm}(x, y) \star {\rm agm}(x, y)$ is also equal to $\theta^2(-q)/\theta^2(q)$ by definition (5),
which completes the proof.  \\ 
\\
\indent
From (B) and (C) of the theorem we see that a positive $\mu$ satisfying $x \star y = \mu \star \mu$ is only
${\rm agm}(x, y)$. Furthermore, (4) and (C) lead to
\[
   \theta^2(q) \star \theta^2(-q)= 1 \star 1 = 1,  
\]
which reveals that  $\theta^2(-q)$ is the inverse of  $\theta^2(q)$ with respect to $\star$.   \\  
\\
\indent
The next theorem states several algebraic properties of the operation, where all variables are assumed to be
positive real numbers. \\
\\
\noindent
{\bf Theorem 2.}  {\it The operation $\star$ satisfies the following algebraic properties.}\\
~{\bf (D)}  $a \star x = a \star y$ {\it implies} $x = y$ ({\it a cancellation law}). \\
~{\bf (E)}  $(ax) \star (ay) = a \star (a(x \star y))$ {\it for any} $a, x, y$ ({\it a distributive law}).\\
~{\bf (F)}  {\it If $z = x \star y$, then} $y  = x(x^{-1} \star (x^{-1}z))$. {\it In particular, the inverse of $x$ 
with respect to the operation is} $x(x^{-1} \star x^{-1})$.\\
\\
\noindent
{\bf Proof.} \\
~(D) Let $\mu = {\rm agm}(a, x)$ and $\mu' = {\rm agm}(a, y)$. If $a \star x = a \star y$, then by (C) we have
$a \star x = \mu \star \mu = a \star y = \mu' \star \mu'$ and hence $\mu = \mu'$ from (B), i.e.,
${\rm agm}(a, x) = {\rm agm}(a, y)$.  This implies ${\rm agm}(1, x/a) = {\rm agm}(1, y/a)$ and hence
$1 \star x/a = 1 \star y/a$ from (C). From (A) we conclude $x/a = y/ a$ and $x = y$. \\
~(E) In the course of the proof of (D) we have obtained the property 
that ${\rm agm}(1, x) = {\rm agm}(1, y)$ entails $x = y$.
Using this and (6) we get $(ax) \star (ay) = a \star (a(x \star y))$ from the following equalities:
\begin{eqnarray*}
{\rm agm}(1,ax \star ay) & = & {\rm agm}(ax, ay) = a{\rm agm}(x, y) = a{\rm agm}(1, x \star y) \\
  & = & {\rm agm}(a, a(x \star y)) = {\rm agm}(1, a \star (a(x \star y))).
\end{eqnarray*}
~(F) The equality $x \star y = z = z \star 1$ implies ${\rm agm}(x, y) = {\rm agm}(z, 1)$ from (6) and
${\rm agm}(ax, ay) = {\rm agm}(az, a)$ for any $a > 0$. This yields $ax \star ay = az \star a$ by (C).
Putting  $a = x^{-1}$,  we obtain $1 \star y/x = z/x \star 1/x$ or $y/x = z/x \star 1/x$. 
Therefore, $y = x(zx^{-1} \star x^{-1})$. In particular, 
letting $z = 1$ leads us to the inverse of $x$. This completes the proof. \\
\\
\indent
As is easily seen from the definition of AGM, ${\rm agm} (x, y) = {\rm agm} ((x + y)/2, \sqrt{xy})$ holds. 
Then (C) implies 
\begin{eqnarray}
     x \star y = (x + y)/2 \star \sqrt{xy}.
\end{eqnarray}
Setting $y = 1$ in (7), we have $x \star 1 = x = (x + 1)/2 \star \sqrt{x}$ for all $x>0$. In particular, 
when $x = (2n+1)^2$ for a natural number $n$,
it is possible to obtain the arithmetic relation for natural numbers:
\[
    (2n +1) \star (2n^2 + 2n + 1)= (2n + 1)^2,
\]
such as $3 \star 5 = 9$, $5 \star 13 = 25$ {\it etc}. It seems that 
this is the only case where three integers are involved. 
It indicates how slowly values of the operation grow compared with multiplication,
but more quickly than addition. \\
\indent
In order to numerically calculate $x \star y$, sequences (1) can be effectively utilized, since they rapidly 
converge to ${\rm agm}(x, y)$.
First, calculate $1/{\rm agm}(x, y)$, which is equal to $A = \theta^2(q)$ for some $q$. 
Second, compute
$B$ satisfying ${\rm agm}(A, B) = 1$ by making use of the bisection method, for example. This $B$ is equal to
$\theta^2(-q)$, that is, $B$ is the inverse of $A$ with respect to the operation. 
It follows from (4) and (5) that the value of $x \star y$ is obtained as the quotient $B/A$. 
Computing in several cases, we see that
the associative law does not hold generally under the operation. \\

\begin{center}
{\bf 3. Some remarks}
\end{center}

Having gotten a required operation, it now enables us to define the same binary operation $x \star y$ 
in another way. 
Here we provide an alternative definition using an elliptic integral instead of the theta function.  
It is given by taking a unique positive solution $x \star y$ of the equation \\
\begin{eqnarray}
\int_0^{\pi/2} \frac{d\varphi}{\sqrt{\cos^2\varphi + (x \star y)^2\sin^2\varphi}}
= \int_0^{\pi/2} \frac{d\varphi}{\sqrt{x^2\cos^2\varphi + y^2\sin^2\varphi}} 
\end{eqnarray}
for each pair of positive numbers $x, y$. The reason why this supplies the same operation comes from the relationship
\[
 \int_0^{\pi/2} \frac{d\varphi}{\sqrt{x^2\cos^2\varphi + y^2\sin^2\varphi}} = \frac{\pi}2 \frac{1}{{\rm agm}(x, y)},
\]
between the elliptic integral and AGM, which was also proved by Gauss. In view of this relationship
we see that (8) is equivalent to the relation ${\rm agm} (1, x \star y) = {\rm agm}(x, y)$
that was obtained from definition (5). Therefore, (8) provides another definition for the same operation $x \star y$. \\ 
\indent
It should be noted that some properties of (A)--(F) in Section 2 are proved easily from the second definition (8). 
Property (A), for example, immediately follows. As for (C), letting $\mu = {\rm agm}(x, y)$, 
it follows that
\[
 \int_0^{\pi/2} \frac{d\varphi}{\sqrt{x^2\cos^2\varphi + y^2\sin^2\varphi}} 
   = \frac{\pi}2 \frac{1}{{\rm agm}(x, y)}
  = \frac{\pi}{2\mu} = \int_0^{\pi/2} \frac{d\varphi}{\sqrt{\mu^2\cos^2\varphi + \mu^2\sin^2\varphi}}.
\]
This implies
\[
\int_0^{\pi/2} \frac{d\varphi}{\sqrt{\cos^2\varphi + (x \star y)^2\sin^2\varphi}}
= \int_0^{\pi/2} \frac{d\varphi}{\sqrt{\cos^2\varphi + (\mu \star \mu)^2\sin^2\varphi}},
\]
proving (C): $x \star y = \mu \star \mu$. \\
\indent
Nevertheless, it seems that the first definition given by (5)
is superior to the second, because it is defined directly and it is helpful 
in deducing relationships between the operation and others.\\
\indent
Finally, we suggest the third for the definition based on hypergeometric series. 
Expanding the integrand 
on the left hand side below, we obtain the expression
\begin{eqnarray*}
 \frac{2}{\pi} \int_0^{\pi /2} \frac{d\varphi}{\sqrt{1 - z\sin^2\varphi}} 
 & = & 1 + \left (\frac{1}{2}\right)^2z + \left(\frac{1\cdot 3}{2 \cdot 4}\right)^2z^2 + 
 \left(\frac{1\cdot 3 \cdot 5}{2 \cdot 4 \cdot 6}\right)^2z^3 + \cdots \\
 & = & F\left (\frac{1}{2}, \frac{1}{2}; 1; z \right).
\end{eqnarray*}
For general hypergeometric series $F(a, b; c; z)$ and their properties we refer to [3]. 
Thus it is evident that the elliptic integral can be replaced by the hypergeometric series in definition (8). 
Since we treat the series only in $|z| < 1$, however, a constraint must be imposed on $x, y$ in defining $x \star y$.  \\
\indent
We show that the operation $x \star y$ can be defined in case of $0 < y \le x <1$, for example. In this case we see 
that $0 < x \star y < 1$.
If we had $x \star y \geq 1$, then equality (8) would not hold, because
\[
 \int_0^{\pi/2} \frac{d\varphi}{\sqrt{\cos^2\varphi + (x \star y)^2\sin^2\varphi}}
< \int_0^{\pi/2} \frac{d\varphi}{\sqrt{x^2\cos^2\varphi + y^2\sin^2\varphi}}. 
\]
Rewriting (8) in terms of the hypergeometric series $F(1/2, 1/2; 1; z)$, it follows that
\[
F\left (\frac{1}{2}, \frac{1}{2}; 1; 1-(x \star y)^2 \right) = 
\frac{1}{x}F\left (\frac{1}{2}, \frac{1}{2}; 1; 1-\frac{y^2}{x^2} \right).
\]   
Thus both series converge and $x \star y$ is well defined.\\
\newpage

\begin{center}
\section*{\normalsize References}
\end{center}
\begin{itemize}
\item[{[1]}] J. Borwein and P. Borwein, {\it Pi and the AGM}, John Wiley, New York, 1987.
\item[{[2]}]  D. A. Cox, The arithmetic-geometric mean of Gauss, {\it L' Enseign. Math.} 30:
275--330, 1984. 
\item[{[3]}]  G. Gasper and M. Rahman, {\it Basic Hypergeometric Series}, Cambridge Univ. Press, Cambridge, 1990.
\item[{[4]}]  S. Tanimoto, A binary operation related to arithmetic-geometric means (Japanese), 
{\it Sugaku} 49: 300--301, 1997. 
\end{itemize}
\end{document}